%% file: article.tex
\documentclass[12pt]{article}
\usepackage[utf8]{inputenc}
\usepackage{xcolor}

\usepackage{natbib} 
\usepackage{ulem}
\usepackage{amsmath,amsthm,amsfonts,dsfont}
\usepackage{mathrsfs}  
\usepackage{tikz}
\usepackage{amssymb}
\usepackage{fullpage}
\usepackage{tabularx}
\usepackage[hidelinks]{hyperref}
\usepackage{float}
\usepackage{subfigure}
\newtheorem{theorem}{Theorem}[section]
\newtheorem{example}{Example}[section]
\newtheorem{property}[theorem]{Property}
\newtheorem{properties}[theorem]{Properties}

\theoremstyle{definition}
\newtheorem{definition}{Definition}[section]
\theoremstyle{remark}
\newtheorem{remark}{Remark}

\usepackage[inline]{enumitem}

\newcommand{\range}[1]{\mathscr{#1}}

\newcommand{\dominantU}{\nu}
\newcommand{\dominantUbar}{\bar{\nu}}

\newcommand{\intensity}{\lambda}
\newcommand{\Intensity}{\Lambda}
\newcommand{\dominantY}{\eta}
\newcommand{\dominantYbar}{\bar{\eta}}
\newcommand{\intensityratio}{\rho}
\newcommand{\densityratio}{\rho}
\newcommand{\parampop}{\theta}
\newcommand{\permutation}{\tau}
\newcommand{\paramnuisance}{\xi}
\newcommand{\provar}{\Sigma}

\newcommand{\subsetA}{A}
\newcommand{\Covariogram}{C}

\newcommand{\Cov}{\mathrm{Cov}}
\newcommand{\derive}{\mathrm{d}}
\newcommand{\Design}{D}
\newcommand{\design}{\mathbf{d}}
\newcommand{\E}{\mathrm{E}}
\newcommand{\density}{\mathrm{f}}
\newcommand{\Semivariogram}{G}

\newcommand{\nrep}{J}
\newcommand{\repindex}{j}
\newcommand{\Sampleindex}{K}
\newcommand{\sampleindex}{\mathbf{k}}
\newcommand{\likelihood}{\mathscr{L}}
\newcommand{\uple}{m}
\newcommand{\Samplesize}{N}
\newcommand{\samplesize}{\mathbf{n}}
\newcommand{\Sample}{S}

\newcommand{\size}{\mathrm{size}}
\newcommand{\Pop}{\mathrm{U}}
\newcommand{\toPop}{\bar{\mathrm{U}}}
\newcommand{\Var}{\mathrm{Var}}

\newcommand{\position}{\mathbf{x}}
\newcommand{\SignalSpace}{\mathscr{Y}}
\newcommand{\toSignalSpace}{\bar{\SignalSpace}}
\newcommand{\Signal}{Y}
\newcommand{\signal}{\mathbf{y}}
\newcommand{\Desvar}{Z}
\newcommand{\DesvarSpace}{\mathscr{Z}}
\newcommand{\desvar}{\mathbf{z}}

\newlength{\defbaselineskip}
\newcommand{\setlinespacing}[1]%
           {\setlength{\baselineskip}{#1 \defbaselineskip}}

 \newcommand{\singlespacing}{\setlength{\baselineskip}{1.3 \defbaselineskip}}

\title{The effect of Informative Selection on the estimation of parameters related to Spatial Processes}
\date{}
\author{Daniel Bonn\'ery \thanks{Epidemiology and Modelling Group, Department of Plant Sciences, University of Cambridge, UK} \and Francesco Pantalone\thanks{Department of Economics, University of Perugia, Italy} \and M. Giovanna Ranalli\thanks{Department of Political Science, University of Perugia, Italy}}

\begin{document}
\normalem 

\maketitle


\input{abstract.tex}
\singlespacing

\input{section1}

\input{section2}
\input{section3}

\input{section4}
\input{sectionC}

\input{Acknowledgements}
\appendix

\input{sectionA.tex}

\bibliographystyle{apalike}
\bibliography{refs}




\end{document}

%% file: abstract.tex
\begin{abstract}
This paper extends the concept of informative selection, population distribution and sample distribution to a spatial process context.
These notions were first defined 
in a context where the output of the random process of interest consists of independent and identically distributed realisations for each individual of a population. 
It has been showed that informative selection was inducing a stochastic dependence among realisations on the selected units. In the context of spatial processes, the ``population'' is a continuous space and realisations for two different elements of the population are not independent. 
We show how informative selection may induce a different dependence among selected units and how the sample distribution differs from the population distribution. 

~

\end{abstract}

%% file: section1.tex
\section{Introduction} \label{sec:intro}
Spatial processes are employed in many fields, such as geology, Earth science, environmental and agricultural surveys, among many others. A huge variety of data is surveyed, such as rainfall data \citep{ord1979spatial}, atmospheric data \citep{thiebaux1987spatial}, forestry data \citep{samra1989spatial} and soils data \citep{burgess1980optimal}, in order to perform statistical inference on parameters of interest. In these applications, a general approach is to postulate a spatial process for the variable $\Signal$ under investigation, that is a stochastic process is assumed to have generated the population values, also called \emph{superpopulation} model. 
When using sample data, we need to pay attention at the relation between the spatial process and the  selection mechanism. Indeed, the distribution of the observed data can be different from the distribution assumed for the population by means of the spatial process $P^{\Signal}$. In other words, the density obtained by the sample data may not be the density obtained by reducing the number of terms in $P^{\Signal}$, as if the units were selected completely at random. In this case, we say that the sampling is \emph{informative}. Ignoring an informative sample may lead to biases and erroneous inference, as illustrated for example in \cite{skinner1989analysis}. This could happen when there is dependence between the assumed stochastic model and the sampling mechanism, that is the units are selected dependently to the variable of interest. Examples of this type of mechanism selection are, among many others, length-biased sampling, endogenous stratification, adaptive sampling, sequential quota sampling and cut-off sampling. For more details, see Section 3 of \cite{bonnery2012uniform}.\\
In this paper, we study the effect of informative selection when a spatial process is employed. We consider informative selection as a situation where the sample responses, \emph{given that they were selected}, are not i.i.d. from the superpopulation model. We start from the notions of population and sample distributions defined by \cite{pfefferman_1992}, and we extend those to a spatial process context. Indeed, \cite{pfefferman_1992} consider realizations of a random process of interest as independent and identically distributed, whereas we address the dependence between the realizations, which very likely characterizes a spatial process. As in \cite{pfeffermann1998parametric}, we start from the distribution of the observed values given they were selected, or \emph{sample pdf}. While \cite{pfeffermann1998parametric} consider the observations as if they were independently distributed according to the sample pdf, our work considers the dependence between the observed units in order to account for the spatial component of the population. Also, while \cite{pfeffermann1998parametric} assume that the target population $U$ is finite, our work can be applied to both finite population and continuous population.\\
Throughout the paper we assume to have a spatial Gaussian process and a point process that represents the selection mechanism. Given these two processes, we then define the probability density function of the sample and a subsample, and the ``sample'' and ``population'' distribution of the signal. These concepts are needed in order to introduce a \emph{density ratio} that we use for investigate the bias eventually introduced by the selection mechanism, or \emph{selection bias}. In doing that, the spatial structure of the population is taken into account.\\
These general definitions allow us to define sample counterpart of different characteristics of the population distribution. Indeed, we focus on the \emph{variogram}, which analyzes the degree of spatial dependence of spatial processes and provides useful insights on the phenomenon under investigation. For instance, the variogram is the key component for the well-known \emph{kriging} method \citep{matheron1962traite}, and it plays a crucial role on prediction since it is used to compute the kriging weights. We introduce the definitions of \emph{population variogram} and \emph{sample variogram}, and we investigate the behaviour of the \emph{naive} estimator, i.e. the estimator that does not take into account the informativeness of the selection mechanism.
\\The paper is organized as follow. In Section \ref{sec:stat_fra} we define the statistical framework used throughout the paper. In particular, general notions and the concepts of spatial process, sample, design and design variables are introduced. In Section \ref{sec:sampledistribution} the sample distribution and the population distribution are defined, and in Section \ref{sec:estimation} estimation of the variogram is briefly reviewed and a small simulation is carried out in order to investigate the behaviour of the naive estimator. Finally, Section \ref{sec:conclusions} provides conclusion and future research.

%% file: section2.tex
\section{Statistical Framework} \label{sec:stat_fra}
This section introduces a statistical framework that describes the relationship between the design and sample on one side, and the variable of interest on the other sides. The following concepts borrow to  the general statistical framework by  \cite{dbb1}.

The definitions and most properties we provide apply in this framework. We also provide a case study example throughout the paper to illustrate the general definitions and properties. Some results presented here only apply to the specific statistical framework used for the example.

\subsection{General notations}

In this paper, all random variables are defined on a probability space $(\Omega,\mathscr{A},P)$. The expected value and variance / covariance operators are defined with respect to the probability measure $P$. 
For two sets $E$ and $F$, $(E\to F)$ or $F^E$ designate the set of the functions from $E$ to $F$.
The notation "$g:E\to F,x\mapsto g(x)$" means: let $g$ be a mapping from set $E$ to set $F$ that to $x$ associates $g(x)$. For $f:E\to (F\to G)$, $g:E\to F$, $h:E\to (H\to F)$,  the notation $f[g]$ designates the function :  $f[g]:E\to G$, such that $f[g](x)=(f(x))(g(x))$, and the notation 
$f[h]$ designates the function :  $f[h]:E\to (H\to G)$, such that $f[h](x)=(f(x)\circ h(x))$.

For a set $E$, $\bar{E}$ designates the set $\{\mathbf{0}\}\cup\bigcup_{n\in\mathbb{N},n\geq 1} E^{\{1,\ldots,n\}}$, where $\{\mathbf{0}\}=E^\emptyset$ corresponds to the set containing the application $\mathbf{0}$ with empty domain. The application $\mathbf{0}$ can be interpreted as the empty sample that matches no draw to the population of interest. 
Let denote by $\size $ the application that maps an application to the cardinality of its domain: $\bar{E}\to\mathbb{N}$, $\mathbf{e}\mapsto n$ if $\position\in E^{\{1,\ldots,n\}}$, $0$ if $\position=\mathbf{0}$. For an application $\position\in\bar{E}$, a set $\Sampleindex$, $\position_\Sampleindex$ is the application:
$\position_\Sampleindex:\Sampleindex\cap\mathrm{domain}(\position)\to E:\ell\mapsto\position(\ell)$, in the case of a random application $\Sample:\Omega\to\toPop$ and a random set :$\Sampleindex\to(\mathscr{P}(\mathbb{N})$ ($\mathscr{P}(\mathbb{N})$ is the set of all subsets of $\mathbb{N}$), then $\Sample_\Sampleindex$ is the random application: $\Omega\to\toPop, 
\omega\mapsto 
\Sample_\Sampleindex(\omega):(\Sampleindex(\omega)\cap\mathrm{domain}(\Sample(\omega))\to \Pop),
\ell\mapsto (\Sample(\omega))(\ell)$.
For a measure $\eta$ of $E$, a non random finite set $\Sampleindex$, $\eta^{\otimes\Sampleindex}$ is the measure such that for any collection $(\subsetA_\ell)_{\ell\in\Sampleindex}$ of subsets of $E$  :
$\eta^{\otimes\Sampleindex}\left(\bigcap_{\ell\in\mathbb{L}}\{\position\in E^\Sampleindex:\position(\ell)\in\subsetA_\ell\}\right)=\prod_{\ell\in\Sampleindex}\eta(\subsetA_\ell)$, and
$\eta^{\otimes\emptyset}(\{\mathbf{0}\})=1$, 
then define the measure $\bar{\eta}$ on $\bar{E}$: $\bar{\eta}=\eta^{\otimes\emptyset}+\sum_{\samplesize\in\mathbb{N},\samplesize\geq 0}\eta^{\otimes\{1,\ldots,\samplesize\}}$.\\
We apply the following rule: Random variables are capital Roman letters, whereas Random variable realisations are bold corresponding lowercase roman letter. By convention, any sum over an empty set is 0 and any product over an empty set is 1.The notation $\density_{V\mid W}$ denotes the density of $V$ conditional on $W$ with respect to a dominating measure on the domain of $V$.

\subsection{Spatial process}
We consider a space $\Pop$, that is a compact 
subset of a finite dimensional real vector space $\mathbb{R}^d$, with its associated Borel sigma-field, and a random process $\Signal$ defined on $\Pop$ with value in another finite dimension real vector space $\SignalSpace$, e.g. $Y:\Omega\to(\Pop\to\SignalSpace)$.
 For example for a random variable $\Sample:\Omega\to\Pop^{\{1,2\}}$, and a random variable
$\Signal:\Omega\to(\Pop\to\SignalSpace)$, 
$\Signal[\Sample]$ is the random variable: $\Signal[\Sample]:\Omega\to(\{1,2\}\to\SignalSpace)$, $\omega\mapsto(\Signal(\omega))(\Sample(\omega)):\ell\to(\Signal(\omega))(\Sample(\omega)(\ell)$. Examples and illustrations will be given for the particular case where $\Pop=[0,1]^2$.
Let $\dominantY$ denote a sigma-finite measure on the set $\SignalSpace$  and let $\dominantU$ denote a probability measure on  $\Pop$.
The {\em $\dominantU$-averaged theoretical semivariogram} is defined as the function: 
\begin{equation}\Semivariogram:\mathbb{R}^d\to[0,+\infty),h\mapsto\frac12\int_{\Pop^{\{1,2\}}} \Var\left[\Signal[\position(2)]-\Signal[\position(1)]\right] \derive(\dominantU^{\otimes \{1,2\}})^{X\mid X[2]-X[1]=h}(\position),\label{eq:averagesemivariogram}\end{equation} 
where $(X)$ is the identity of $\Pop^{\{1,2\}}$, and the {\em $\dominantU$-averaged theoretical covariogram} is the function  $\nu^{X_2-X_1}-a.s(h)$-defined:
\begin{equation}\Covariogram:\mathbb{R}^d\to\mathbb{R}, h\mapsto\int_{\Pop^{\{1,2\}}} \Cov\left[\Signal[\position(1)],\Signal[\position(2)]\right] \derive(\dominantU^{\otimes \{1,2\}})^{X\mid X[2]-X[1]=h}(\position).\label{eq:averagecovariogram}\end{equation} The covariogram and semivariogram satisfy the relationship: $\forall h\in\mathbb{R}^d, \Semivariogram(h)=\Covariogram(0)-\Covariogram(h)$.

\begin{definition}[Intrinsic stationarity and Second order stationarity, \protect{\citep[p.~53]{cressie2015statistics}}]
A process is {\em intrinsic stationary} when the following conditions are satisfied :  $\forall \position\in\Pop^{\{1,2\}}$,
\begin{eqnarray}
    \mathrm{E}\left[\Signal\left[\position(2)\right]-\Signal\left[\position(1)\right]\right]&=&0\\
    \frac12~\Var\left[\Signal\left[\position(2)\right]-\Signal\left[\position(1)\right]\right]&=&\Semivariogram(\position(2)-\position(1))\label{eq:semivariogram}
\end{eqnarray}
A process is {\em second order stationary} when the following are satisfied:
\begin{equation}
\exists\mu\in\mathbb{R},~    \forall\position\in \Pop,~~ E\left[\Signal\left[\position\right]\right]=\mu
\end{equation}
\begin{equation} \label{eq:covariogram}
    \forall~\position\in \Pop^{\{1,2\}}, ~\Cov\left[\Signal\left[\position(1)\right],\Signal\left[\position(2)\right]\right]=C\left(\position(1)-\position(2)\right)
\end{equation}
\end{definition}
In the case of a first order stationary process, the variance operator $\Var[.]$ can equivalently be replaced by the square expected value operator $\mathrm{E}[(.)^2]$ in equation \eqref{eq:averagesemivariogram}.
In the case of a second order stationary process, equations \eqref{eq:averagecovariogram} and \eqref{eq:averagesemivariogram} correspond to the definition of the theoretical covariogram and semivariogram as found in \citet[p.~53 and p.~58]{cressie2015statistics}.
The random process is isotropic if in addition on being second order stationary, the covariogram function $h\mapsto \Covariogram(h)$ only depends on $h$ via $h\mapsto\|h\|$.


Common model assumptions on the process $\Signal$ consist in assuming second order stationarity and isotropy. Covariance structure of the signal is then fully characterized by $\Covariogram(0)$ and $\Semivariogram(h), h\neq 0$. 
For the simulations and illustrations in this paper, will be used
a Gaussian covariogram (see \citet[p.~80]{chiles1999geostatistics} for more models) of the form
\begin{equation}\label{def:gaussiancovariogram}\Covariogram(h)=\parampop_1\exp\left(-\|h\|^2/\parampop_2^2)\right)\end{equation} where the deviation parameter $\parampop_1$, and the scale parameter $\parampop_2$ are real positive numbers.




\begin{definition}[Isotropic Gaussian Random Process]
\label{def:isotropicgaussianprocess}
A {\em Gaussian random process} $\Signal:\Omega\to(\Pop\to\SignalSpace)$ is such that its distribution  can be derived from the distributions of $\Signal[\position]$, where $\position\in\toPop$. 
For $\position$, $\position'\in\toPop$, denote the expected value of the signal by $\mu:\toPop\to\toSignalSpace,\position\mapsto\mathrm{E}\left[\Signal[\position]\right]$, and the covariance of the random vectors $\Signal[\position]$, $\Signal[\position']$  by $\provar_{\Signal;\position,\position'}=\Cov \left[\Signal[\position],\Signal[\position']\right]$.
The distribution of $\Signal$ is then fully characterized by $\mu$ and $\provar$: for $n\in\mathbb{N}$,   $\position\in\Pop^{\{1,\ldots,\samplesize\}}$, $\Signal[\position]$ has the following density with respect to $\dominantY^{\otimes \{1,\ldots,n\}}$: 
\begin{equation} \label{eq:pdf_norm_process}
    \density_{\Signal[\position]}\left(\signal\right)=\left(2\pi^{n/2}|\provar_{\Signal;\position,\position}|^{\frac12}\right)^{-1}\exp\left(-\frac12(\signal-\mu(\position))\provar_{\Signal;\position,\position}^{-1}(\signal-\mu(\position))^{\!\mathrm{T}}\right).
\end{equation}
The distribution of an  {\em isotropic Gaussian process} is such that $\forall \position, \position'\in\Pop$,  $\provar_{\Signal;\position,\position'}= \Covariogram\left(\|\position-\position'\|\right)$ and is fully characterized by $\mu$ and $\Covariogram$.
\end{definition}

\begin{example}[Isotropic Gaussian process $\Signal$ with Gaussian Covariogram]\label{example:main}

We simulate three independent replications of an isotropic Gaussian process with $\Signal:\Omega\to(\Pop=[0,1]^2\to\mathbb{R})$, with $\forall \position\in\Pop$, $\mathrm{\mu}[\position]=0$, with a Gaussian Covariogram \eqref{def:gaussiancovariogram} with deviation parameter
 $\parampop_1=1$, and scale parameter $\parampop_2=\lbrace 0.01, 0.1, 1\rbrace$, respectively.
Figure \ref{fig:oaijsfdoij} represents independent realizations of $\Signal$.
\begin{figure}[H]
    \caption{Heat maps of realisations of the random process $\Signal:\Omega\to(\Pop=[0,1]^2\to\mathbb{R})$ for different values of $\parampop_2$}
    \label{fig:oaijsfdoij}
    
\hspace{-.6cm}\includegraphics{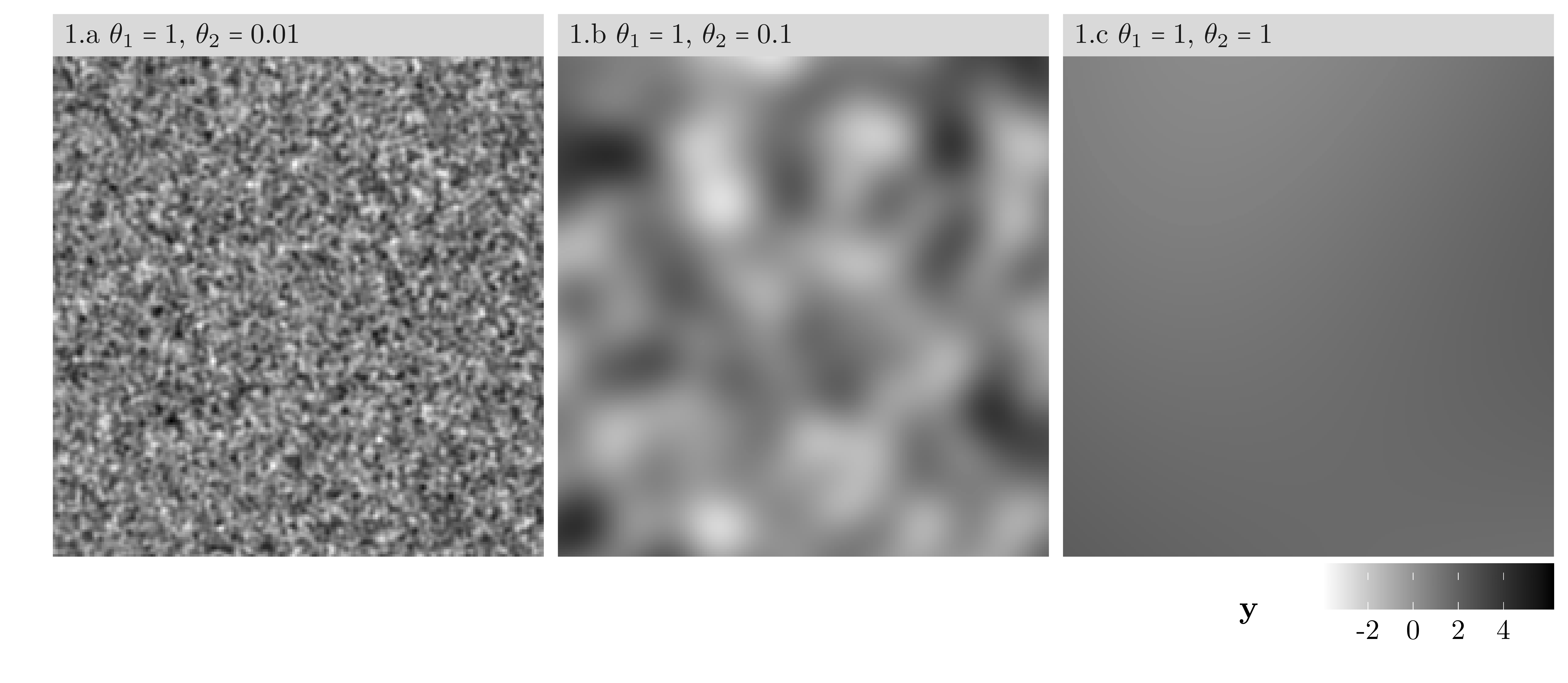}

    \vspace{-.4cm}
    {\footnotesize For each realisation $\signal$ of the random process $\Signal$, for each $\position$ in $\Pop$, the value of $\signal[\position]$ is color coded with a grayscale gradient.}
\end{figure}

\end{example}

\subsection{Sample, design and design variable} \label{sec:design}

\subsubsection{Fixed designs and design variables}

By definition a design $\design$ is a probability distribution on $\toPop$.
A sample $\Sample$ drawn from $\design$ is a random variable or point process of distribution $\design$, e.g. a random variable $\Sample$ such that $\mathrm{P}^\Sample=\design$. 
Define the size $\Samplesize=\mathrm{size}\circ\Sample$ of the sample $\Sample$. The sample density with respect to $\dominantUbar$ is defined by:
$(\derive P^\Sample)/(\derive\dominantUbar)(\position)=P(N=n)\times(\derive P^{\Sample\mid N=n})/(\derive\dominantU^{\otimes\{1,\ldots,n\}})(\position)$, if $\position\in\Pop^{\{1,\ldots,n\}}$, $P(N=0)$ if $\position=\mathbf{0}$. 
A fixed design variable is a function $\desvar:\Pop\mapsto \DesvarSpace$.
A fixed size design $\design$ is usually defined as a function of a fixed design variable, and characterised by its density with respect to $\dominantUbar$. For example, the Probability Proportional to Size $\desvar$ With Replacement and size $n$ design, with $\DesvarSpace=[0,+\infty)$, for $\position\in \toPop$, is characterized by:  
\begin{equation}\label{eq:ppswr}\density_{\Sample}(\position)=\left(\derive P^{\Sample}/\derive\dominantUbar\right)(\position)=\left|\begin{array}{l}\left(\int_U \left(\desvar(\position')\right)\mathrm{d}\dominantU (\position')\right)^{-\samplesize}
\left(\prod_{\ell=1}^\samplesize\left(\desvar(\position(\ell))\right)\right)\text{ if }\position\in\Pop^{\{1,\ldots,n\}},\\~0\text{ otherwise.}
\end{array}\right.\end{equation}

We remind that with our notations, $(\desvar.\dominantU)(\Pop)=\int_\Pop \desvar(\position)\derive \dominantU(\position)$.The point process $\Sample$ characterized by Equation \eqref{eq:ppswr} is a binomial point process of $\samplesize$ points in $\Pop$ with intensity $\Pop\to\mathbb{R}, \position\to\left(\left(\desvar.\dominantU\right)\left(\Pop\right)\right)^{-1}\desvar(\position)$, which we abbreviate by $\Sample\sim\mathrm{bpp}\left(\desvar,\samplesize\right)$. Simple random sampling with replacement is a binomial point process with a constant intensity,  and we refer to as $\Sample\sim\mathrm{bpp}\left(1,\samplesize\right)$.  
For example, given a measurable real function $\desvar :\Pop \to\mathbb{R}$, a spatial Poisson Process of intensity $\desvar$ is a point process 
$S:\Omega\to\bigcup_{\samplesize\in\mathbb{N}}\Pop^{\samplesize}$, such that 
for all $\dominantU$-measurable subset $\subsetA$ of $\Pop$, \begin{equation}\label{eq:poissonprocess}
S\sim\mathrm{Ppp}(\desvar)\Leftrightarrow
\mathrm{cardinality}(\Sample^{-1}[\subsetA])\sim \mathrm{Poisson}\left((\desvar.\dominantU)\left(\subsetA \right)\right),
\end{equation}

where $\Sample^{-1}[\subsetA]$ is the random variable with domain the finite subsets of $\Pop$ defined by: $\omega\mapsto\{\ell\in\{1,\ldots,\Samplesize(\omega)\};(\Sample(\omega))(\ell)\in \subsetA\}$ if $\Samplesize(\omega)>0$, $\emptyset$ otherwise.

The density of such process with respected to $\bar\dominantU$ is defined, for $\position\in \toPop$, by:
\begin{equation}
\density_{\Sample}(\position)=\left(\derive P^{\Sample}/\derive\dominantUbar\right)(\position)=\left(\size (\position)!\right)^{-1}\exp\left(-\left(\desvar.\dominantU\right)(\Pop)\right)
\prod_{\ell\in\mathrm{domain}(\position)}
\desvar(\position(\ell)).
\end{equation}

For a point process $\Sample$ with values in the measured space $(\Pop,\dominantU)$, for a random variable $W$, define $\Intensity_{S\mid W=w}$ as the intensity measure of $\Sample$ conditionally to $W=w$ with respect to the measure on $\Pop$, where for each measurable subset $\subsetA$ of $\Pop$, $\Intensity_S(A)=\mathrm{E[\mathrm{cardinality}(S^{-1}[A])\mid W=w]}$, and $\intensity_S$ as the density of $\Intensity_S$ with respect to $\dominantU$: $\intensity_S=\derive\Intensity_S/\derive\dominantU$.

\subsubsection{Random design variables, sample and random design}
In practice, the design parameter $\desvar$ is modeled as the output of a random process $\Desvar:\Omega\to(\Pop\to \DesvarSpace)$ that we will refer to as the design variable.
The selection process, when controlled, is in practice a function of an auxiliary variable, called design variable, that is a process defined on the same space $\Pop$. When the selection process is not chosen by the experimenter it can also be modelled as a function of such a process that can be observed, partially observed or latent. In practice, it may not be reasonable to assume independence of the design variable  $\Desvar$ and study variable $\Signal$.

The design is by definition a random variable with domain the set of probability distributions on  $\toPop$. The sample is a random variable $\Sample$ with domain $\toPop$ such that the distribution of $\Sample$ conditionnaly to the design \emph{is} the design, e.g:
\begin{equation}P^\Design-a.s.(\design),~P^{\Sample\mid \Design=\design}=\design.\end{equation}

We assume that the distribution of the sample conditionally on the design, the signal and the design variable is the design:
\begin{equation}\label{eq:SconditionalonDindependentonYZ}
P^{\Design,\Signal,\Desvar}-a.s.(\design,\signal,\desvar),~P^{\Sample\mid \Design=\design,\Signal=\signal,\Desvar=\desvar}=\design.\end{equation}
This assumption is an independence of the sample $\Sample$ on the design variable and the signal conditionally on the design, which does not imply independence of the sample and the signal.



%

 For simplicity, we only consider exchangeable index sample designs, in the sense that $\forall \samplesize\in\mathbb{N}$, for all permutation $\permutation$ of $\{1,\ldots,n\}$,
 \begin{equation}P^\Desvar-\text{a.s.}(\desvar), P^{S\mid\Samplesize=\samplesize,\Desvar=\desvar}=P^{S[\permutation]\mid \Samplesize=\samplesize,\Desvar=\desvar}\label{assumption:designindexexchangeability}.\end{equation}
which ensures that 
 \begin{equation}P^{S\mid\Samplesize=\samplesize}=P^{S[\permutation]\mid \Samplesize=\samplesize}\label{assumption:designindexexchangeability}.\end{equation}

\begin{example}[Example \ref{example:main} continued: distribution of $\Desvar$ conditionally on $\Signal$]
\label{example:2.2}
We assume that 
$\Desvar$ satisfies:
$$\Desvar=\mathrm{exp}\left(\paramnuisance_0+\paramnuisance_1\Signal+\paramnuisance_2 \varepsilon\right),$$
where $\paramnuisance_0$,  $\paramnuisance_1$, $\paramnuisance_2$ are real positive numbers, 
$\varepsilon:\Omega\to\left(\Pop\to\mathbb{R}\right)$
is an isotropic Gaussian process with mean
$\mu=0$ and Gaussian covariogram $\Covariogram(h)=\paramnuisance_{\text{deviation}}\exp\left(-\|h\|^2/\paramnuisance_{\text{scale}}^2)\right)$ with deviance parameter   
 $\paramnuisance_{\text{deviation}}=1$ and scale parameter $\paramnuisance_{\text{scale}}$. 

The variable $\signal:\Pop\to\mathbb{R}$ (resp $e:\Pop\to\mathbb{R}$) is generated once by sampling from $\Signal$ (resp $\varepsilon$). The variable $\desvar=\exp\left(\paramnuisance_0+\paramnuisance_1\signal+\paramnuisance_2\mathbf{e}\right)$ is computed for 3 different values of the vector 
$(\paramnuisance_0,\paramnuisance_1, \paramnuisance_2)$, chosen so that the expected value of the sampling intensity is the same infor each simulation, more specifically so that $E[\log(\Desvar)]=exp(\paramnuisance_0+((\parampop_1\times\paramnuisance_1)^2+\paramnuisance_2^2)/2)=10$. Two samples are drawn. The first sample, conditionally on $\Desvar=\desvar$, follows $P^{\Sample\mid\Desvar=\desvar}=\mathrm{bpp}(\desvar,10)$. The unconditional distribution of $\Sample$  is denoted $P^\Sample=\mathrm{bpp}(\Desvar,10)$ in this case. The second sample, conditionally on $\Desvar=\desvar$, follows $P^{\Sample\mid\Desvar=\desvar}=\mathrm{Ppp}(\desvar)$. The unconditional distribution of $\Sample$  is denoted  $P^\Sample=\mathrm{Ppp}(\Desvar)$.
The variable $\desvar$, generated with different parameters is mapped in Figure 2.a, 2.b, 2.c, and the variable $\signal$ is mapped in figure 2.d. In particular, Figures 2.a, 2.b, 2.c. show how sampled units tend to concentrate where the sampling intensity $\desvar$ is the highest.
By comparing Figure 2.b and Figure 2.d, as $\paramnuisance_1$ is 0, the high sampling intensity in Figure 2.b and sampled points does not necessarily correspond to high values of the signal $\signal$ in Figure 2.d. 
In the opposite, the high sampling intensity in Figure 2.c and sampled points often corresponds to high values of the signal $\signal$ in Figure 2.d. 

\begin{figure}[H]
    \caption{Heat maps of the design variable $\desvar$ and plot of realisations of $\Sample$ for three different design variables and designs.}
    \label{fig:oaijsfdwefweoij}
    \hspace{-.6cm}\includegraphics{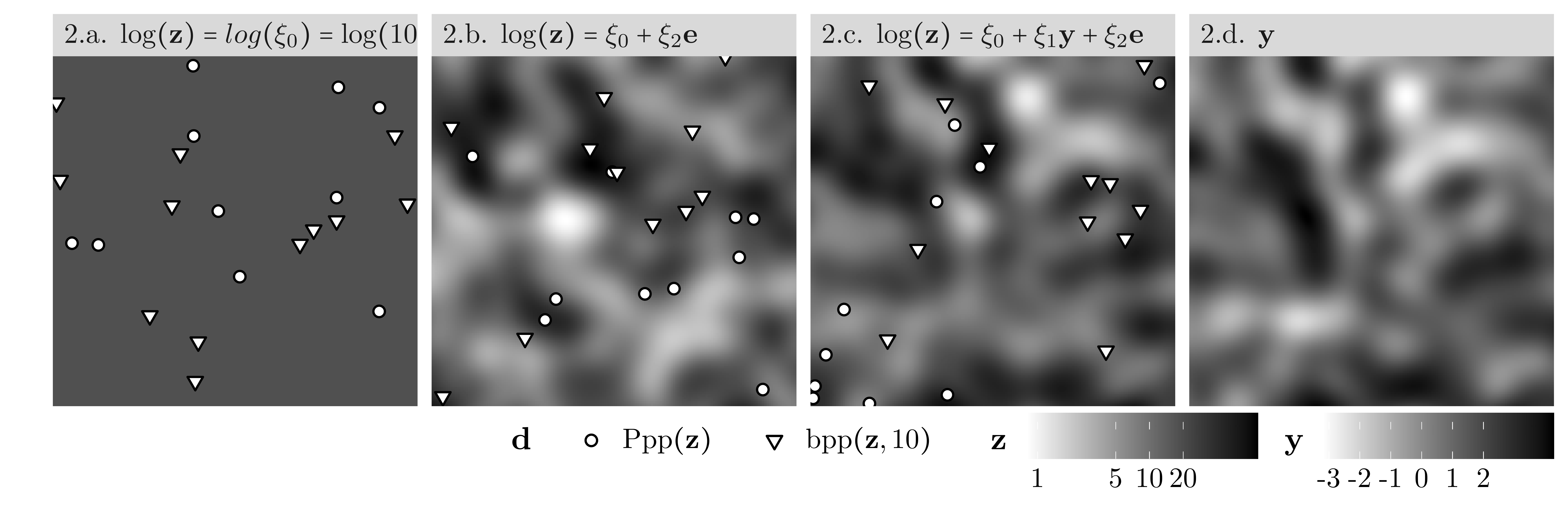}
    \vspace{-1cm}

    \footnotesize{

Sub-figures correspond to the  heatmap of a design variable $\desvar$ (Sub-figures 2.a, 2.b, 2.c) to the heatmap of $\signal$ (Sub-figure 2.d) and the plot of the samples (circle and triangle dots) drawn according to the two different designs. The values of $(\paramnuisance_0, \paramnuisance_1 ,\paramnuisance_2)$ for each sub-figure are: 2.a: $(\log(10),0,0)$, 2.b:$(\log(10)-0.125,0,0.5)$, 2.c: $(\log(10)-0.125,0.4,0.3)$.}

\end{figure}

Figure \ref{fig:jointmarginaldensitiesYZ} represents the values of $\signal(\position)$ vs $\desvar(\position)$ for values of $\position$ on a regular grid of $\Pop$.
The values of $\desvar$ in Figure \ref{fig:jointmarginaldensitiesYZ}.a (resp 
\ref{fig:jointmarginaldensitiesYZ}.b, 
\ref{fig:jointmarginaldensitiesYZ}.c) correspond to the values of Figure 2.a (resp 
Figure 2.b, 
Figure 2.c). Figure \ref{fig:jointmarginaldensitiesYZ} illustrates that high sampling intensity correspond to high values of $\signal$ when $\paramnuisance_2>0$, as in Figure  \ref{fig:jointmarginaldensitiesYZ}.c.

\begin{figure}[H]
\caption{Joint and marginal densities of $\desvar$ and $\signal$}\label{fig:jointmarginaldensitiesYZ}
\includegraphics{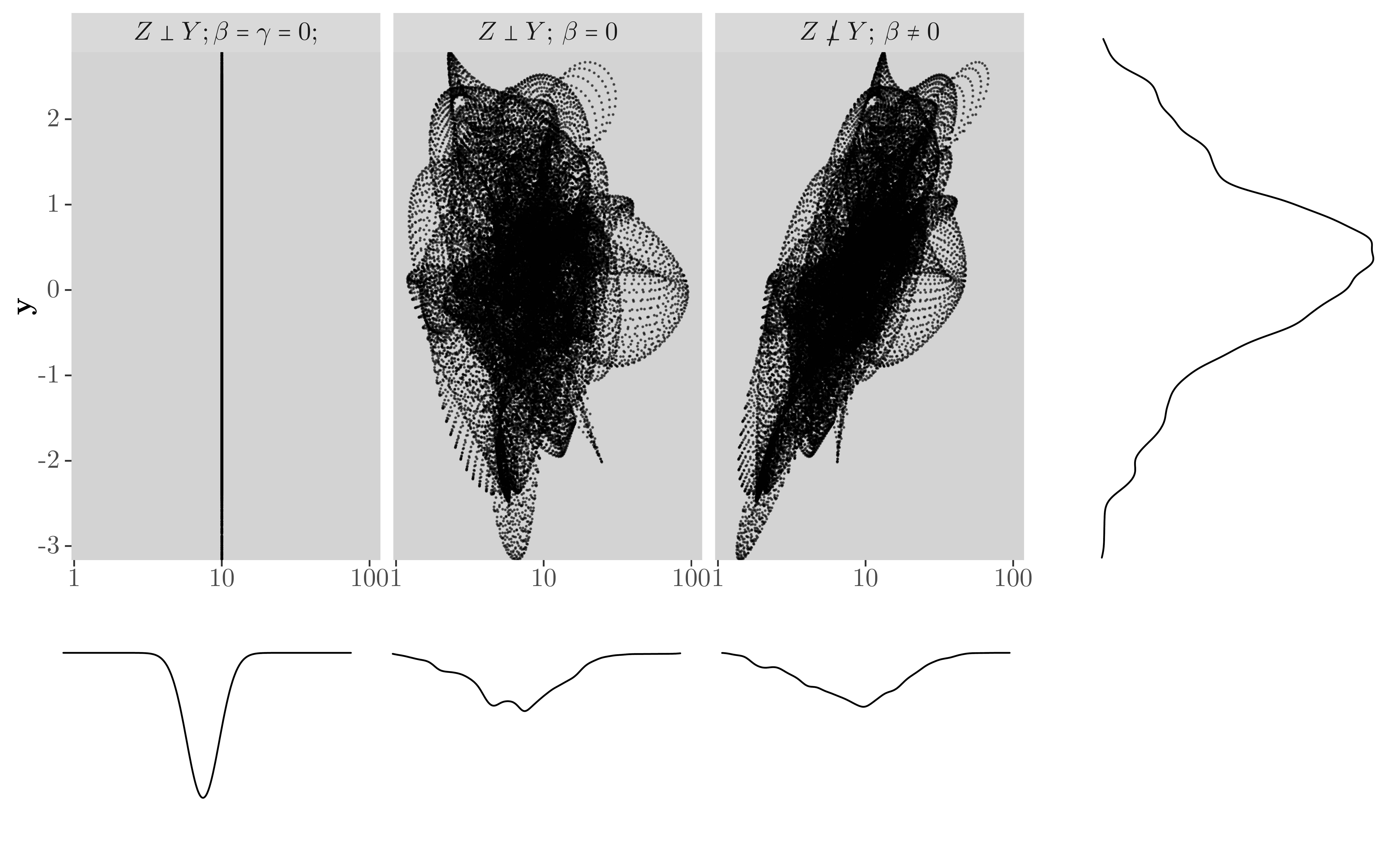}
\footnotesize{
 Each sub-figure contains the scatter plot of  $\left\{\left(\desvar(\position),\signal(\position)\right):\position\in\mathrm{Grid}\right\}$, where $\mathrm{Grid}$ is a regularly spaced grid of $\Pop$. The values of $(\paramnuisance_0, \paramnuisance_1 ,\paramnuisance_2)$ for each Sub-figure are: 3.1: $(\log(10),0,0)$, 3.2:$(\log(10)-(0.5^2+0.3^2),0,\sqrt{0.5^2+0.3^2})$, 3.3: $(\log(10)-(0.5^2+0.3^2),0.5,0.3)$. The vertical axis corresponds to $\signal$, the vertical to $\desvar$. The marginal plots correspond to the density of $\desvar$ (right margin) and to the densities of $\signal$ (bottom margins).}
\end{figure}
\end{example}


\subsection{Probability density function of the sample}

The conditional distributions of the sample and of subsamples can be derived from  the conditional distribution of the design:

\begin{eqnarray}
\density_{\Sample\mid W}(\position\mid\mathbf{w})&=&
\int \left(\derive P^{\Sample\mid \Design,W}/\derive\dominantUbar \right)(\position\mid \design,\mathbf{w})\derive P^{\Design\mid W}(\design\mid\mathbf{w})
\end{eqnarray}

In particular when $W$ is a function of $\Signal$ and $\Desvar$, assumption \eqref{eq:SconditionalonDindependentonYZ} implies that the probability to select the sample $\position$ conditionally on $W=\mathbf{w}$ is the 
probability to select the sample $\position$ given the design, averaged on all possible designs conditionally on  $W=\mathbf{w}$:
\begin{eqnarray}
\density_{\Sample\mid W}(\position\mid\mathbf{w})&=&
\int \left(\derive \design/\derive\dominantUbar \right)(\position)~\derive P^{\Design\mid W}(\design\mid\mathbf{w}).
\end{eqnarray}.

\begin{example}[ Example \ref{example:2.2} continued: Sample density for $\Sample\sim \mathrm{Ppp}(\Desvar)$ and $\Sample\sim \mathrm{bpp}(\Desvar,n)$]
\label{example:2.3}
 For $\position\in\Pop^{\{1,\ldots,\samplesize\}}$, when $\Sample\sim \mathrm{Ppp}(\Desvar)$,
\begin{equation}\density_{\Sample\mid W}(\position|\mathbf{w})=
\int\exp\left(-(\desvar.\dominantU)(\Pop)\right)(\samplesize!)^{-1}\left(\prod_{\ell=1}^\samplesize \desvar(\position(\ell))\right)\derive
P^{\Desvar\mid W}(\desvar\mid \mathbf{w}).\label{eq:oijoijejods}
\end{equation}
 When $\Sample\sim \mathrm{bpp}(\Desvar,\samplesize)$,
\begin{equation}
\density_{\Sample\mid W}(\position|\mathbf{w})=
\int\left(-(\desvar.\dominantU)(\Pop)\right)^{-\samplesize}\left(\prod_{\ell=1}^\samplesize \desvar(\position(\ell))\right)\derive
P^{\Desvar\mid W}(\desvar\mid \mathbf{w}).\label{eq:oijerogiql}
\end{equation}
\end{example}

\subsection{Probability density function of a subsample}

Theoretical developments in Section \ref{sec:sampledistribution} require to express the density function of a random sized or fixed size subsample of a random or fixed size sample.

For a fixed set $\Sampleindex$, the density 
$\density_{\Sample_\Sampleindex\mid W}(\position\mid\mathbf{w})$
can be derived from $\density_{\Sample\mid W}$ via:
\begin{eqnarray}
\density_{\Sample_\Sampleindex\mid W}(\position\mid\mathbf{w})&=&
\int \density_{\Sample\mid W}(\position'\mid \mathbf{w})~\mathds{\mathds{1}}(\position'_\Sampleindex==\position)~\derive\dominantUbar(\position').
\end{eqnarray}

More generally, for a random set $\Sampleindex$ the density 
$\density_{\Sample\mid W}(\position\mid\mathbf{w})$
can be derived from $\density_{\Sample_\Sampleindex\mid W}$ via:
\begin{eqnarray}
\density_{\Sample_\Sampleindex\mid W}(\position\mid\mathbf{w})&=&
\int\left(\int \density_{\Sample\mid W}(\position'\mid \mathbf{w})~\mathds{\mathds{1}}(\position'_\Sampleindex==\position)~\derive\dominantUbar(\position')\right)\derive P^{K\mid W=\mathbf{w}}.
\end{eqnarray}

\begin{remark}\label{remark:K1subsetK2}
For  two fixed size sets  $\Sampleindex_1, \Sampleindex_2$ fixed such that 
$\Sampleindex_1\subset \Sampleindex_2$ we do not have necessarily, the equation 
\begin{equation}\label{eq:iugheirguhei}\density_{\Sample_{\Sampleindex_1}}(\position\mid \mathbf{w})=\int \density_{\Sample_{\Sampleindex_2}}(\position'\mid \mathbf{w})\mathds{1}(\position'_{\Sampleindex_1}=\position)\derive \dominantU^{\otimes K_2}(\position')
\end{equation}
is not generally true. 
For example, when $S\sim \mathrm{Ppp}(\Desvar)$, 
$\Sampleindex_1=\{1\}$, $\Sampleindex_2=\{1,2\}$,
For $\position \in \Pop^{\{1\}}$, 
$\density_{\Sample_{\{1\}}}(\position)=\density_{\Sample}(\position)+\int \density_{\Sample_{\Sampleindex_2}}(\position'\mid \mathbf{w})\mathds{1}(\position'_{\Sampleindex_1}=\position)\derive \dominantU^{\otimes K_2}(\position')$.
A sufficient condition for the equation \eqref{eq:iugheirguhei} would be that $\Samplesize$ is not random.
\end{remark}
\begin{example}[Example \ref{example:2.3} continued: Probability density function of a subsample for $\Sample\sim \mathrm{Ppp}(\Desvar)$ and $\Sample\sim \mathrm{bpp}(\Desvar,\samplesize)$ .]
\label{example:2.4}

Let $\Sampleindex$ be  a finite subset of $\mathbb{N}\setminus\{0\}$. Let   $\Sampleindex'$ be an element of $\{\{\ell\in \Sampleindex:\ell\leq\samplesize\}:\samplesize\in\mathbb{N}\}$, and let  $\position\in\Pop^{\Sampleindex'}$. Let $W$ be a function of $\Desvar$, $\Signal$.
In the following we derive the density of $\Sample_\Sampleindex$ conditionally on $W=\mathbf{w}$.
For $\position=\mathbf{0}$,

\begin{equation}\density_{\Sample_K\mid W}(\mathbf{0}\mid\mathbf{w})=\left|\begin{array}{ll}1&\text{if }\Sampleindex=\emptyset,\\
P(\Samplesize<\min(K)\mid W=\mathbf{w})&\text{ otherwise.}
\end{array}\right.
\end{equation}
For $\position\neq\mathbf{0}$:
\begin{equation}
\density_{\Sample_K\mid W}(\position\mid\mathbf{w})=
\int \density_{\Sample_K\mid \Desvar}(\position\mid\desvar) \derive P^{\Desvar\mid W}(\desvar\mid\mathbf{w}).
\end{equation}

When $\Sample\sim \mathrm{Ppp}(\Desvar)$,
\begin{equation}
\density_{\Sample_K\mid \Desvar}(\position\mid\mathbf{\desvar})=
\left(\sum_{\substack{\samplesize\in\mathbb{N}:\\\{\ell\in\Sampleindex:\ell\leq\samplesize\}=\Sampleindex'}}\!\!\!\!\!\!\!\! P(\Samplesize=\samplesize\mid\Desvar=\desvar)\right)\prod_{\ell\in\Sampleindex'}\frac{\desvar(\position(\ell))}{(\desvar.\dominantU)(\Pop)}.
\end{equation}

and when $\Sample\sim \mathrm{bpp}(\Desvar,\samplesize)$,
\begin{equation}
\density_{\Sample_{\Sampleindex}\mid W}(\position|\mathbf{w})=\left|\begin{array}{ll}
\int\left(\prod_{\ell\in\Sampleindex'} \left(\left((\desvar.\dominantU)(\Pop)\right)^{-1}\left(\desvar(\position(\ell))\right)\right)\right)\derive
P^{\Desvar\mid W}(\desvar\mid \mathbf{w})&\text{if }\Sampleindex'=\{1,\ldots,\samplesize\}\cap\Sampleindex\\
0&\text{otherwise.}
\end{array}\right.\end{equation}

\end{example}

Note that the distribution of the sample unconditionally on the design or the design variable can be obtained via: 
$\density_{\Sample}(\position)=\int \density_{\Sample\mid\Desvar}(\position)\derive P^\Desvar$.
Contrary to the particular case of population index exchangeability and sample index exchangeability, the sample index exchangeability \eqref{assumption:designindexexchangeability} condition alone does not ensure that conditionally on the sample size, $\Samplesize=\samplesize$, the distribution of the sample  $P^{\Sample\mid\Samplesize=\samplesize}$ is the uniform distribution on $\Pop^{\{1,\ldots,\samplesize\}}$. 
For example, $f_{\Sample_{\{1\}}}$ can be uniform, but $f_{\Sample_{\{1,2\}}}$ not.

%% file: section3.tex
\section{``Sample" vs ``Population''  distribution of the signal}\label{sec:sampledistribution}

The observation consists of the realisations of the random variables $\Sample$ and $\Signal[\Sample]$. In this section, we derive the distribution of the observed values of the signal on the sample, (e.g. the distribution of $\Signal[\Sample]$) from  the distribution of the design variable conditionally $\Desvar$ to the signal $\Signal$ and the function that links the design to the design variable, or equivalently the distribution of the sample $\Sample$ conditionally to the design variable $\Desvar$. We resort to  the Bayes formula to express the density of $(\Signal[\Sample],\Sample)$ in $(\signal,\position)$,   as the product of the density of the sample $\Sample$ in $\position$ conditionally on $\Signal[\position]=\signal$ by the density of the signal $\Signal[\position]$ in $\signal$. 
\subsection{Density ratio}
\subsubsection{Definition and properties}
We start with the following definition.


\begin{definition}
For a random set $\Sampleindex$, 
define
$\densityratio_{\Sampleindex}(.\mid.)$ as any function that satisfies:

\begin{equation}
P^{(\Sample_\Sampleindex,\Signal[\Sample_\Sampleindex])}-\text{a.s}(\position,\signal),~\density_{\Signal[\position]\mid\Sample_\Sampleindex=\position}\left(\signal\right)=
    \density_{\Signal[\position]}\left(\signal\right)
    \densityratio_{\Sampleindex}\left(\position\mid  \signal\right)\label{eq:owijoj}
\end{equation}

\end{definition}

\begin{property}\label{prop:ieuhierhgi}
\begin{equation}
P^{(\Sample_\Sampleindex,\Signal[\Sample_\Sampleindex])}-\text{a.s}(\position,\signal),~~~
\densityratio_{\Sampleindex}\left(\position \mid \signal\right)
         \density_{\Sample_\Sampleindex}\left(\position\right)=
    \density_{\Sample_\Sampleindex\mid \Signal[\position]}\left(\position|\signal\right)
\end{equation}

\end{property}

\begin{proof}
From the Bayes formula:
\begin{eqnarray}
\density_{\Signal[\position]\mid\Sample_\Sampleindex =\position }\left(\signal\right)
&=&(\density_{\Sample_\Sampleindex }(\position))^{-1}~\density_{\Signal[\position],\Sample_\Sampleindex }(\signal,\position)\\
&=&(\density_{\Sample_\Sampleindex }(\position))^{-1}~
\density_{\Sample_\Sampleindex \mid\Signal[\position]}(\position\mid \signal)~
\density_{\Signal[\position]}(\signal)\label{eq:oijoijsf}
\end{eqnarray}

So, combining Equations \eqref{eq:oijoijsf} and  \eqref{eq:owijoj}:
\begin{equation}
{\densityratio_{\Sampleindex}\left(\position \mid \signal\right)}
=(\density_{\Sample_\Sampleindex }(\position))^{-1}
~\density_{\Sample_\Sampleindex \mid\Signal[\position]}(\position\mid \signal)
\end{equation}
\end{proof}%

Property \ref{prop:ieuhierhgi} provides an interpretation of $\densityratio_\Sampleindex(\signal\mid \position)$ as the probability that the draws indexed by $\Sampleindex$ correspond to the sample $\position$ conditionally on the values of $\Signal$ in $\position$ divided by the same probability unconditionally on the values of $\Signal$. 
In this sense, the density ratio $\densityratio_\Sampleindex(\signal\mid \position)$ is a relative probability of selection.

\begin{property}
The density ratio $\densityratio$ can be derived from the distribution of $\Desvar$ via:

\begin{eqnarray}\label{eq:rho_expression}
\densityratio_{\Sampleindex}(\position\mid\signal)&=&\left({\int
         \density_{\Sample_\Sampleindex\mid \Desvar}\left(\position|\desvar\right)\derive P^{\Desvar}}\right)^{-1}{\int\density_{\Sample_\Sampleindex\mid \Desvar}(\position\mid\desvar)\derive P^{\Desvar\mid\Signal[\position]=\signal}(\desvar)}
\label{eq:oigjeogjio}\end{eqnarray}
\end{property}

\begin{property}\label{prop:indeprho}
If $\Desvar$ and $\Signal$ are independent, then $\densityratio_K(\position\mid\signal)=1$.
\end{property}
\begin{proof}
$\Desvar\perp\Signal\Rightarrow P^\Desvar=P^{\Desvar\mid\Signal[\position]=\signal}$. Replacing $P^{\Desvar\mid\Signal[\position]=\signal}$ by $P^{\Desvar}$ in the numerator of Equation \eqref{eq:rho_expression} yields  $\densityratio_K(\position\mid\signal)=1$.
\end{proof}

\subsubsection{Computation of $\densityratio$ on a Case study}
In this subsection, we use the framework of Example \ref{example:2.4} to illustrate the properties and provide graphical representations of the density ratio.
We remind that in this framework, $\Sample\sim \mathrm{Ppp}(\Desvar)$ or  $\Sample\sim \mathrm{bpp}(\Desvar,\samplesize)$.
Conditionally on $\Signal[\position]=\signal$, $\Desvar$ is a log normal random process with distribution characterized by 
$E[\log(\Desvar[\position'])\mid \Signal[\position]=\signal]=\paramnuisance_0+\paramnuisance_1(\mu+\Sigma_{\position',\position}\Sigma_{\position,\position}^{-1} (\signal-\mu))$ and 
$\mathrm{Var}[\log(\Desvar[\position'])\mid \Signal[\position]=\signal]=\paramnuisance_2^2\provar_{\varepsilon;\position',\position'}+\paramnuisance_1^2\left(\provar_{\Signal;\position',\position'}-\provar_{\Signal;\position',\position}\provar_{\Signal;\position,\position}^{-1}\provar_{\Signal;\position,\position'})\right)$.
\begin{proof}
See Appendix \ref{sec:A.owifjoeij}
\end{proof}

For $\position=\mathbf{0}$  and  $\signal=\mathbf{0}$,
\begin{eqnarray*}
\lefteqn{\density_{\Sample_{\{1\}}\mid\Signal[\mathbf{0}]=\mathbf{0}}(\mathbf{0}))}\\&=&
\density_{\Sample_{\{1\}}}(\mathbf{0})\\
&=&\left|\begin{array}{ll}
\int \exp\left(-(\desvar.\dominantU)(\Pop)\right)\derive P^{\Desvar}(\desvar)&\text{if }\Sample\sim\mathrm{Ppp}(\Desvar)\\
0&\text{if }\Sample\sim\mathrm{bpp}(\Desvar,\samplesize),\end{array}\right.
\end{eqnarray*}
so $\densityratio_{\{1\}}(\mathbf{0}\mid\mathbf{0})=
\density_{\Sample_{\{1\}}\mid\Signal[\mathbf{0}]=\mathbf{0}}(\mathbf{0})/\density_{\Sample_{\{1\}}}(\mathbf{0})=
1$.
More generally, for any random set $\Sampleindex$, $\densityratio_{\Sampleindex}(\mathbf{0}\mid\mathbf{0})=1.$
When $\Sample\sim \mathrm{Ppp}(\Desvar)$,
for $\position\in\Pop^{\{1\}}$, 
\begin{eqnarray*}
\lefteqn{\densityratio_{\{1\}}(\position\mid\signal)}\\&=&\frac{
\int 
\left(\sum_{\samplesize\geq 1} (\samplesize!)^{-1}(\desvar.\dominantU)(\Pop)^\samplesize\right)\times \exp\left(-(\desvar.\dominantU)(\Pop)\right)
\left((\desvar.\dominantU)(\Pop)\right)^{-1}\desvar(\position(1))\derive P^{\Desvar\mid \Signal[\position]}(\desvar\mid\signal) 
}{
\int 
\left(\sum_{\samplesize\geq 1} (\samplesize!)^{-1}(\desvar.\dominantU)(\Pop)^\samplesize\right)\times \exp\left(-(\desvar.\dominantU)(\Pop)\right)
\left((\desvar.\dominantU)(\Pop)\right)^{-1}\desvar(\position(1))\derive P^{\Desvar}(\desvar) 
}\\&=&\frac{
\int 
\left(1-\exp\left(-(\desvar.\dominantU)(\Pop)\right)\right)
\left((\desvar.\dominantU)(\Pop)\right)^{-1}\desvar(\position(1))\derive P^{\Desvar\mid \Signal[\position]}(\desvar\mid\signal) 
}{
\int 
\left(1-\exp\left(-(\desvar.\dominantU)(\Pop)\right)\right)
\left((\desvar.\dominantU)(\Pop)\right)^{-1}\desvar(\position(1))\derive P^{\Desvar}(\desvar)\phantom{11111} 
}\\&=&\frac{
\int 
\left(1-\exp^{-\int \exp\left( \paramnuisance_0+\paramnuisance_1\Signal[\position']+\paramnuisance_2\varepsilon[\position']\right)\derive\dominantU(\position')}\right)
\left(\int \exp^{-\paramnuisance_1(\Signal[\position']-\signal)-\paramnuisance_2(\varepsilon[\position']-\varepsilon[\position])}\derive\nu(\position')\right)^{-1}\derive P^{\Signal,\varepsilon\mid \Signal[\position]=\signal} 
}{
\int
\left(1-\exp^{-\int \exp\left( \paramnuisance_0+\paramnuisance_1\Signal[\position']+\paramnuisance_2\varepsilon[\position']\right)\derive\dominantU(\position')}\right)
\left(\int \exp^{-\paramnuisance_1(\Signal[\position']-\Signal[\position(1)])-\paramnuisance_2(\varepsilon[\position']-\varepsilon[\position(1)])}\derive\nu(\position')\right)^{-1}\derive P^{\Signal,\varepsilon}\phantom{11111} 
}
\end{eqnarray*}

and for  $\position\in\Pop^{\{1,2\}}$, 
\begin{eqnarray*}
\lefteqn{\densityratio_{\{1,2\}}(\position\mid\signal)}\\&=&\frac{
\int 
\frac{\left(1-(1+(\int\exp\left( \paramnuisance_0+\paramnuisance_1\Signal[\position']+\paramnuisance_2\varepsilon[\position']\right)\derive\dominantU(\position')))\right)}{\exp\left(\int\exp\left( \paramnuisance_0+\paramnuisance_1\Signal[\position']+\paramnuisance_2\varepsilon[\position']\right)\derive\dominantU(\position')\right)}\frac{
\exp^{-\paramnuisance_1(\sum_{\ell=1}^2\signal(\ell))-\paramnuisance_2(\sum_{\ell=1}^2\varepsilon[\position(\ell)])}}{
\left(\int \exp\left(-\paramnuisance_1(\Signal[\position'])-\paramnuisance_2(\varepsilon[\position'])\right)\derive\nu(\position')\right)^2
}\derive P^{\Signal,\varepsilon\mid \Signal[\position]=\signal} 
}{
\int 
\frac{\left(1-(1+(\int\exp\left( \paramnuisance_0+\paramnuisance_1\Signal[\position']+\paramnuisance_2\varepsilon[\position']\right)\derive\dominantU(\position')))\right)}{\exp\left(\int \exp\left(\paramnuisance_0+\paramnuisance_1\Signal[\position']+\paramnuisance_2\varepsilon[\position']\right)\derive\dominantU(\position')\right)}\frac{
\exp^{-\paramnuisance_1(\sum_{\ell=1}^2\Signal(\position[\ell]))-\paramnuisance_2(\sum_{\ell=1}^2\varepsilon[\position(\ell)])}}{
\left(\int \exp\left(-\paramnuisance_1(\Signal[\position'])-\paramnuisance_2(\varepsilon[\position'])\right)\derive\nu(\position')\right)^2
}\derive P^{\Signal,\varepsilon}\phantom{11111} 
}
\end{eqnarray*}

and for  $\position\in\Pop^{\{1,\ldots,\samplesize\}}$, 
\begin{eqnarray*}
\lefteqn{\densityratio_{\{1,\ldots,\Samplesize\}}(\position\mid\signal)}\\&=&\frac{
\int 
\frac{\left(\int \exp\left(\paramnuisance_0+\paramnuisance_1\Signal[\position']+\paramnuisance_2\varepsilon[\position']\right)\derive\dominantU(\position')\right)^\samplesize}{\samplesize!~\exp\left(\int \paramnuisance_0+\paramnuisance_1\Signal[\position']+\paramnuisance_2\varepsilon[\position']\derive\dominantU(\position')\right)}\frac{
\exp^{-\paramnuisance_1(\sum_{\ell=1}^\samplesize\signal(\ell))-\paramnuisance_2(\sum_{\ell=1}^\samplesize\varepsilon[\position(\ell)])}}{
\left(\int \exp\left(-\paramnuisance_1(\Signal[\position'])-\paramnuisance_2(\varepsilon[\position'])\right)\derive\nu(\position')\right)^\samplesize
}\derive P^{\Signal,\varepsilon\mid \Signal[\position]=\signal} 
}{
\int 
\frac{\left(\int\exp\left( \paramnuisance_0+\paramnuisance_1\Signal[\position']+\paramnuisance_2\varepsilon[\position']\right)\derive\dominantU(\position')))\right)^\samplesize}{\samplesize!\exp\left(\int \exp\left(\paramnuisance_0+\paramnuisance_1\Signal[\position']+\paramnuisance_2\varepsilon[\position']\right)\derive\dominantU(\position')\right)}\frac{
\exp^{-\paramnuisance_1(\sum_{\ell=1}^2\Signal(\position[\ell]))-\paramnuisance_2(\sum_{\ell=1}^\samplesize\varepsilon[\position(\ell)])}}{
\left(\int \exp\left(-\paramnuisance_1(\Signal[\position'])-\paramnuisance_2(\varepsilon[\position'])\right)\derive\nu(\position')\right)^\samplesize
}\derive P^{\Signal,\varepsilon}\phantom{11111} 
}
\end{eqnarray*}

When $\Sample\sim \mathrm{bpp}(\Desvar,\samplesize)$, for $\position \in \Pop^{\{1\}}$:

\begin{equation*}
\densityratio_{\{1\}}(\position\mid\signal)=
\frac{
\int 
\left(\int \exp^{-\paramnuisance_1(\Signal[\position']-\signal)-\paramnuisance_2(\varepsilon[\position']-\varepsilon[\position])}\derive\nu(\position')\right)^{-1}\derive P^{\Signal,\varepsilon\mid \Signal[\position]=\signal} 
}{
\int
\left(\int \exp^{-\paramnuisance_1(\Signal[\position']-\Signal[\position(1)])-\paramnuisance_2(\varepsilon[\position']-\varepsilon[\position(1)])}\derive\nu(\position')\right)^{-1}\derive P^{\Signal,\varepsilon}\phantom{11111} 
}.
\end{equation*}

for $\position \in \Pop^{\{1,2\}}$:
\begin{equation*}
\densityratio_{\{1,2\}}(\position\mid\signal)=
\frac{
\int 
\frac{\exp\left(-\paramnuisance_1(\sum_{\ell=1}^2\signal(\ell)-\paramnuisance_2(\sum_{\ell=1}^2\varepsilon[\position(\ell)])\right)}{\left(\int \exp\left(-\paramnuisance_1\Signal[\position']-\paramnuisance_2\varepsilon[\position']\right)\derive\nu(\position')\right)^{2}}\derive P^{\Signal,\varepsilon\mid \Signal[\position]=\signal} 
}{
\int 
\frac{\exp\left(-\paramnuisance_1(\sum_{\ell=1}^2\Signal[\position(\ell)]-\paramnuisance_2(\sum_{\ell=1}^2\varepsilon[\position(\ell)])\right)}{\left(\int \exp\left(-\paramnuisance_1\Signal[\position']-\paramnuisance_2\varepsilon[\position']\right)\derive\nu(\position')\right)^{2}}\derive P^{\Signal,\varepsilon}\phantom{11111} 
}.
\end{equation*}

and 
for $\position \in \Pop^{\{1,\ldots,\samplesize'\}}$, $\samplesize'\leq\samplesize$:

\begin{equation*}
\densityratio_{\{1,\ldots,\samplesize'\}}(\position\mid\signal)=
\frac{
\int 
\frac{\exp\left(-\paramnuisance_1(\sum_{\ell=1}^2\signal(\ell)-\paramnuisance_2(\sum_{\ell=1}^{\samplesize'}\varepsilon[\position(\ell)])\right)}{\left(\int \exp\left(-\paramnuisance_1\Signal[\position']-\paramnuisance_2\varepsilon[\position']\right)\derive\nu(\position')\right)^{\samplesize'}}\derive P^{\Signal,\varepsilon\mid \Signal[\position]=\signal} 
}{
\int 
\frac{\exp\left(-\paramnuisance_1(\sum_{\ell=1}^{\samplesize'}\Signal[\position(\ell)]-\paramnuisance_2(\sum_{\ell=1}^{\samplesize'}\varepsilon[\position(\ell)])\right)}{\left(\int \exp\left(-\paramnuisance_1\Signal[\position']-\paramnuisance_2\varepsilon[\position']\right)\derive\nu(\position')\right)^{\samplesize'}}\derive P^{\Signal,\varepsilon}\phantom{11111} 
}.
\end{equation*}

The density ratio is a function of $\position$, $\signal$,  $\paramnuisance$, and $\parampop$.
\begin{properties}{Properties of $\densityratio$ common to Poisson and binomial point processes}
\begin{enumerate}
\item[i.] $\paramnuisance_1=0\Rightarrow \intensityratio(\position,\signal;\paramnuisance,\parampop)=1$.
\item[ii.] $\intensityratio_\Sampleindex(\position\mid\signal;\paramnuisance,\parampop)-1=o_{\paramnuisance_2\to\infty}(1)$,
\item[iii.] $\intensityratio_\Sampleindex(\position\mid\signal;\paramnuisance,\parampop)-1=o_{\parampop_{\text{scale}}\to\infty}(1)$,
\item[iv.] $\intensityratio_\Sampleindex(\position\mid\signal;\paramnuisance,\parampop)-1=o_{\parampop_{\text{scale}}\to\infty}(1)$,
\item[v.] $\intensityratio_\Sampleindex(\position,-\signal;(\parampop,(\paramnuisance_0,-\paramnuisance_1,\paramnuisance_2,\paramnuisance_{\text{deviation}},\paramnuisance_{\text{scale}})=\intensityratio_\Sampleindex(\position\mid\signal;(\parampop,(\paramnuisance_0,\paramnuisance_1,\paramnuisance_2,\paramnuisance_{\text{deviation}},\paramnuisance_{\text{scale}})$
\item[vi.] $\mathrm{sign}(\paramnuisance_1)=1\Rightarrow \lim_{\signal\to\infty}\intensityratio_{\{1\}}(\position\mid\signal;\paramnuisance,\parampop)=+\infty$,
\item[vii.] $\mathrm{sign}(\paramnuisance_1)=1\Rightarrow \lim_{\signal\to-\infty}\intensityratio_{\{1\}}(\position\mid\signal;\paramnuisance,\parampop)=0$,
\end{enumerate}
\end{properties}

\begin{properties}{Properties of $\densityratio$ when $\Sample\sim\mathrm{bpp}(\Desvar,\samplesize)$}
\begin{enumerate}
\item[i.] $ \intensityratio_{\{1\}}(\position\mid\signal;\paramnuisance,\parampop)=
\exp(\paramnuisance_1(\signal-\paramnuisance/2) + o_{\parampop_{\text{scale}}\to 0}(1)$.
\item[ii.] $\samplesize'\leq\samplesize\Rightarrow\intensityratio_{\{1,\ldots\samplesize'\}}(\position\mid\signal;\paramnuisance,\parampop)=
\exp\left(\paramnuisance_1(\sum_{\ell=1}^{\samplesize'}\signal(\ell)-\samplesize'\paramnuisance_1/2)\right) + o_{\parampop_{\text{scale}}\to 0}(1)$.
\end{enumerate}
\end{properties}




The density ratio $\densityratio(\position\mid\signal)$ is the product of  two integrals as shown in Equation \eqref{eq:oigjeogjio}. No close form could be obtained for any of these two integrals. However we can use Monte Carlo approximations to compute $\densityratio(\position\mid\signal)$.
In the case where $\Pop$ is finite, we can simulate realisations of $P^\Signal$,  and  $P^\varepsilon$ on the population independently $\nrep$ times.
Let denote by $\Signal^{(\repindex)}$  and $\varepsilon^{(\repindex)}$ the $j$-th realisation of each distribution  $\Signal$ and $\varepsilon$ respectively.

For $\repindex \in 1,\ldots,\nrep$, let define $\Signal_c^{(\repindex)}$ by 
$\Signal_c^{(\repindex)}[\position']=\Signal^{(\repindex)}[\position']+\provar_{\Signal;\position',\position}\provar_{\Signal;\position,\position}^{-1}(\signal-\Signal^{(\repindex)}[\position])$, which ensures that 
$P^{\Signal_c^{(\repindex)}}=P^{\Signal\mid\Signal[\position]=\signal}$ by the Matheron rule \citep{wilson2020pathwise, doucet2010note}.

Define, for a non random set $\Sampleindex$,  $\emptyset\neq\Sampleindex\subset\{1,\ldots,\samplesize\}$, 
and $\position\in\Pop^\Sampleindex$, $\signal\in\SignalSpace^\Sampleindex$, $\uple=\size(\Sampleindex)$:
$$\tilde\densityratio_{\Sampleindex}(\position\mid\signal)=
\frac{\sum_{\repindex=1}^\nrep 
    \left(\left(
    \left(\exp\left(\paramnuisance_1\Signal_c^{(\repindex)}+\paramnuisance_2\varepsilon^{(\repindex)}\right).\dominantU\right)(\Pop)\right)^{-\uple}\prod_{\ell\in\Sampleindex}\exp\left(\paramnuisance_1\signal(\ell)+\paramnuisance_2\varepsilon[\position(\ell)]\right)
    \right)}{\sum_{\repindex=1}^\nrep 
    \left(\left(
    \left(\exp\left(\paramnuisance_1\Signal^{(\repindex)}+\paramnuisance_2\varepsilon^{(\repindex)}\right).\dominantU\right)(\Pop)\right)^{-\uple}\prod_{\ell\in\Sampleindex}\exp\left(\paramnuisance_1\Signal_{\phantom{c}}^{(\repindex)}[\position(\ell)]+\paramnuisance_2\varepsilon[\position(\ell)]\right)
    \right)}.$$

In the case of a binomial point process ($\Sample\sim\mathrm{bpp}(\Desvar,\samplesize)$), the ratio of the expected values of the numerator and denominators of $\tilde\densityratio(\position\mid\signal)$ is exactly $\densityratio_\Sampleindex(\position\mid\signal)$.
The empirical variance of the summands of the numerator and numerator can be used to estimate the variance of $\tilde\densityratio(\position\mid \signal)$.

\begin{figure}[H]
\centering
\caption{Representation of $\tilde\densityratio_{\lbrace 1\rbrace}(\position\mid\signal)$ when $\Sample \sim \mathrm{bpp}(\Desvar,\samplesize)$, $\paramnuisance_{2} = 1$, and $\paramnuisance_{1}$ varies. The value of $\signal$ is reported on top of the corresponding graph.}
\includegraphics[width = 0.7\textwidth]{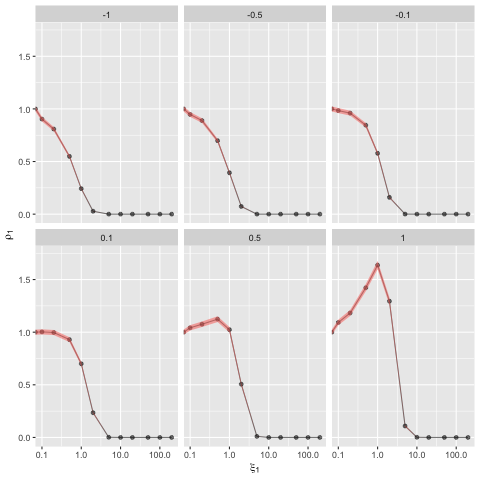}
\label{fig:rho_1_beta}
\end{figure}

\begin{figure}[H]
\centering
\caption{Representation of $\tilde\densityratio_{\lbrace 1\rbrace}(\position\mid\signal)$ when $\Sample \sim \mathrm{bpp}(\Desvar,\samplesize)$, $\paramnuisance_{1} = 1$, and $\paramnuisance_{2}$ varies. The value of $\signal$ is reported on top of the corresponding graph.}
\includegraphics[width = 0.7\textwidth]{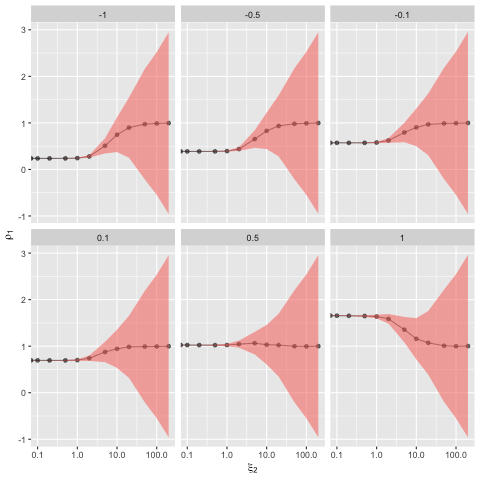}
\label{fig:rho_1_gamma}
\end{figure}

\begin{figure}[H]
\centering
\caption{Representation of $\tilde\densityratio_{\lbrace 1,2\rbrace}(\position\mid\signal)$ when $\Sample \sim \mathrm{bpp}(\Desvar,\samplesize)$, $\paramnuisance_{1} = 1$, $\paramnuisance_{2}=1$, and for four different couples of points.}
\includegraphics[width = 0.8\textwidth]{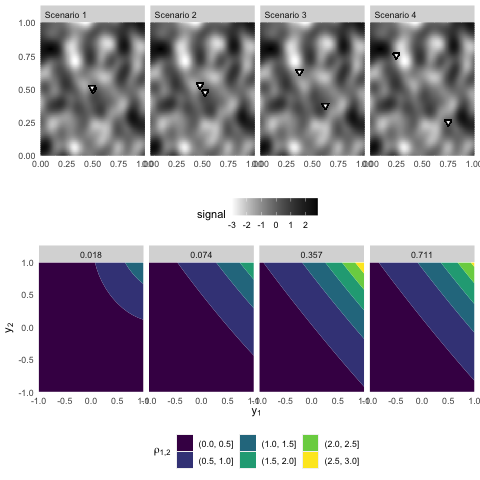}
\label{fig:rho_12}
\end{figure}

Following the aforementioned discussion, we employ a Monte Carlo simulation in order to approximate $\densityratio$ by means of $\tilde\densityratio_{\Sampleindex}(\position\mid\signal)$. As first step, we generate $J=10,000$ replications for the computation of $\tilde\densityratio_{\lbrace 1\rbrace}(\position\mid\signal)$. In Figure \ref{fig:rho_1_beta} we investigate the role of $\paramnuisance_{1}$. Toward this end, $\paramnuisance_{2}$ is kept fixed to 1, and $\paramnuisance_{1}$ varies while $\signal$ is kept fixed as well. This procedure is done for different levels of $\signal$. As expected, the value of $\tilde\densityratio_{\lbrace 1\rbrace}(\position\mid\signal)$ at $\paramnuisance_{1}=0$ is always equal to 1. Indeed, we recall that $\paramnuisance=0$ implies independence between $\Signal$ and $\Desvar$. When $\paramnuisance_{1}$ differs from zero, $\tilde\densityratio_{\lbrace 1\rbrace}(\position\mid\signal)$ is different from 1, hence reflecting the fact that dependence between $\Signal$ and $\Desvar$ introduces bias. In particular, we can observe that as $\paramnuisance_{1}$ approaches infinity, $\tilde\densityratio_{\lbrace1\rbrace}(\position\mid\signal)$ approaches zero. Also, for positive values of $\signal$, the maximum is achieved around $\paramnuisance_{1}=\signal$. The same analysis is performed for the behavior of $\paramnuisance_{2}$, where $\paramnuisance_{1}$ is kept fixed to 1, and different values of $\signal$ are fixed while $\paramnuisance_{2}$ varies. We report the results in Figure \ref{fig:rho_1_gamma}. When $\paramnuisance_{2}=0$, $\tilde\densityratio_{\lbrace 1\rbrace}(\position\mid\signal)$ is different from zero. This reflects the fact that $\paramnuisance_{1}$ is fixed to 1. As $\paramnuisance_{2}$ increases, the effect of $\paramnuisance_{1}$ decreases and $\tilde\densityratio_{\lbrace 1\rbrace}(\position\mid\signal)$ tends to 1. We then turn the attention to $\tilde\densityratio_{\lbrace 1,2\rbrace}(\position\mid\signal)$. $J=4,000$  replications are generated and values of $\tilde\densityratio_{\lbrace 1,2\rbrace}(\position\mid\signal)$ are computed for four different couples of units, which have an Euclidean distance (between the two points) of 0.018, 0.074, 0.357, and 0.711, respectively. In Figure \ref{fig:rho_12} the contours of the surface generated by the values of $\tilde\densityratio_{\lbrace 1,2\rbrace}(\position\mid\signal)$ when the two values of $\signal$ vary are plotted.

\subsection{Distribution of \texorpdfstring{$\Signal[\Sample]$}{Y[\Sample_\Sampleindex]} and weighted density}
The density of $\Signal[\Sample_\Sampleindex]$ with respect to $\dominantYbar$ is defined by
\begin{eqnarray}\density_{\Sample_\Sampleindex,\Signal[\Sample_\Sampleindex]}(\position\mid\signal)&=& \density_{\Signal[\position]}(\signal)~\times~\densityratio_{\Sampleindex}(\position\mid\signal)~\times~\density_{\Sample_\Sampleindex}(\position).\label{eq:iuhsdewripeowi}\\
&=& \density_{\Signal[\position]}(\signal)~\times~\density_{\Sample_\Sampleindex\mid\Signal[\position]}(\position\mid\signal).\label{eq:iuhsdewripeowi2}
\end{eqnarray}

So when $\Desvar\perp\Signal$,
\begin{equation}\density_{\Signal[\Sample_\Sampleindex]\mid\Sample_\Sampleindex}(\signal\mid \position)= \density_{\Signal[\position]}(\signal).\label{eq:iuhsdewripeowi2}\end{equation}

\cite{pfefferman_1992} explore  the case where
$K=\{1\}$ and the spatial process on satisfies the  population independence assumption:

$\forall \position\in\Pop^{\{1,\ldots,\samplesize\}}$, such that $(\position[1],\ldots,\position[\samplesize])$ are two by two  distinct,

\begin{equation}
(\Signal[\position[1]],\ldots,\Signal[\position[\samplesize]])\text{ are i.i.d variables}.
\label{eq:independenceassumption}
\end{equation}

To the 'population' distribution, which is the distribution of $\Signal$ that satisfies Equation \eqref{eq:independenceassumption} with probability density function $\density_{\Signal[\position]}$ for $\position \in\Pop^{\{1\}}$, \cite{pfefferman_1992} oppose the sample distribution, that is the  distribution of a random variable $\Signal^\star$ that satisfies Equation \eqref{eq:independenceassumption}  and with probability density function  $\densityratio_{\{1\}}\density_{\Signal[\position]}$ for $\position \in\Pop^{\{1\}}$.

The random variable $\Signal^\star$ ''does not exist``, in the sense that if the observer could measure $\Signal[\position]$ for all values of $\position$ on $\Pop$, of for values of $\position$ in a non informative sample (e.g. such that $\densityratio=1$), the distribution of  observations would have probability density $\density_\Signal[\Sample]$.
With an informative sample, the observations are similar to what one would observe with a simple random sample and if the population was following $P^{\Signal^\star}$ and not $P^\Signal$. The distribution of $\Signal^\star$ is called by \cite{pfefferman_1992} the ``sample'' distribution.

In the case of a spatial process, can we define a sample distribution ?

Sometimes not, due to the property pointed out in Remark \ref{remark:K1subsetK2}:
the finite dimensional probability density functions 
$\densityratio{\{1,\ldots,\samplesize\}}(\position\mid\signal)\density_{\Signal[\position](\signal)}$ are not necessarily the finite dimensional densities of the same random process distribution.
A sufficient condition for the finite dimensional densities to be finite dimensional probability functions of the same process is that the sample size is fixed.
But then, there is not necessarily unicity of the random process distribution that have these finite dimensional probability density functions, as they are only defined up to the dimension equal to the sample size.

Generalising \cite{pfefferman_1992}, the selection is non informative when $\forall \uple\in\mathbb{N}$ $P^{\Signal[\Sample],\Sample}-\text{a.s.}(\signal,\position)$, $$\densityratio_{\{1,\ldots,\Samplesize\}}(\position\mid\signal)=1.$$ In this case, we consider that the sample distribution corresponds to the population distribution.

In the general case, we then define the sample distribution as a
collection of finite dimensional distributions with respective p.d.f.'s
$\densityratio{\{1,\ldots,\samplesize'\}}(\position\mid\signal)\density_{\Signal[\position](\signal)}$, with\\ $\samplesize'\in\left\{\{1,\ldots,\max\left(\Samplesize(\Samplesize^{-1}(\mathbb{N}))\right)\right\}$.

\subsection{Sample counterpart of distribution characteristics}

Even if the definition of the sample distribution is not the distribution of a random process on $\Pop$, we propose to define the sample counterpart of different characteristics of the population distribution as for example 
the sample semivariogram.

In a non spatial framework,  \cite{pfefferman_1992} differentiates the characteristics of the distribution of one variable when observed on the sample, and the characteristics of the same variable when observed on the population. The terminology to differentiate conists in adding "population" or "sample" to the name of the characteristic. For example, the population probability density function (pdf) is the equivalent  $\density_{\Signal[\position]}$ in our example, and the sample pdf is the density $\density_{\Signal[\Sample[1]]}$. \cite{bonnery2012uniform} and \cite{dbb1} investigate further this concept of population and sample characteristics.
In the following, we provide the definition of the sample characteristics of the distribution of $\Signal$ that are meaningful in the context of a spatial model.

Some characteristics of the distribution of $\Signal$ can be  defined as expected values of a transformation of the values of $\Signal$ on all samples of $\Pop$.
For example, we have seen that the distribution of 
a Gaussian process is characterised by $E[\Signal[\position]]$ for $\position\in \Pop$ and
$\Var[\Signal[\position]]$ for $\position\in \Pop^{\{1,2\}}$. Those characteristics can be obtained by integrals of the form $\int g(\signal[\position]) \density_{\Signal[\position]}(\signal[\position])\derive\dominantYbar$, which themselves can be seen as expected values of statistics derived from the values of $\Signal$ on a simple random sample with replacement.
Those characteristics are called ``population characteristics''. Their sample counterparts are obtained by substituting $ \density_{\Signal[\position]}$ by $\densityratio_{\{1,\ldots,\size(\position)\}}\density_{\Signal[\position]}$, and can be interpreted as the expected value of the same statistics derived from the values of $\Signal$ on the informative sample $\Sample$.





In the case of sampling with fixed size $\samplesize\geq 2$, 
define the sample semivariogram for exchangeable designs as 
$$\Semivariogram^\star(h)=\frac12~\E\left[\left(\Signal[\Sample[1]]-\Signal[\Sample[2]]\right)^2\mid \Sample[2]-\Sample[1]=h \right].$$

\begin{property}[Relationship between $\Semivariogram$, $\Semivariogram^\star$ and $\densityratio$]
$$\Semivariogram(h)^\star=\frac12\int_{\Pop^{\{1,2\}}} \left[\int_{\range{\Signal}^{\{1,2\}}} (\signal(2)-\signal(1))^2~ \density_{\Signal[\position]}(\signal)~
\intensityratio_{\{1,2\}}(\position\mid\signal)~
\derive\dominantY^{\otimes 2}(\signal)
\right] \derive(\dominantU^{\otimes\{1, 2\}})^{X\mid X[2]-X[1]=h}(\position)$$
\end{property}

In the case of a random size sample, one can define the sample covariogram as 
$$\Semivariogram^\star(h)=\frac12\sum_{\samplesize\geq2}P(\Samplesize=\samplesize)~\E\left[\left(\Signal[\Sample[1]]-\Signal[\Sample[2]]\right)^2\mid \Samplesize=\samplesize \text { and }\Sample[2]-\Sample[1]=h \right].$$

\begin{proof}
Let $\position\in\Pop^{\{1,2\}}$, then
\begin{eqnarray*}
\lefteqn{\E\left[\left(\Signal[\position(2)]-\Signal[\position(1)]\right)^2\mid \Sample_{\{1,2\}}=\position\right]}\\
&=&\int_{\range{\Signal}^2} (\signal(2)-\signal(1))^2~
\density_{\Signal[\position]\mid S_{\{1,2\}}}(\signal\mid \position)~
\derive\dominantY^{\otimes\{1, 2\}}(\signal) \\
&=&\int_{\range{\Signal}^2} (\signal(2)-\signal(1))^2~
\densityratio_{\{1,2\}}(\position\mid\signal)~ \density_{\Signal[\position]}(\signal)~
\derive\dominantY^{\otimes\{1, 2\}}(\signal) \\
\end{eqnarray*}
\end{proof}

%
%




%% file: section4.tex
\section{Estimation} \label{sec:estimation}

\subsection{Naive non parametric estimation : the case of the variogram}
With the term \emph{naive estimation} we refer to the situation where the estimator does not take into account the mechanism selection. Therefore, in case of informative sampling, the estimator could suffer from selection bias.
A common approach to achieve a valid variogram estimator is composed by \emph{estimation} and \emph{fitting} of the variogram. In the former, an estimate of the variogram is obtained, while the latter phase is necessary since the estimators used in the first phase are usually not conditionally negative-definite.

Under the assumption of constant-mean, an estimator based on the method of moments is \citep{matheron1962traite}
\begin{equation*} \label{eq:variogram_hat}
2\hat{\gamma}\left(h\right)=\frac{1}{|N\left(h\right)|}\sum_{(\ell_1,\ell_2)\in N\left(h\right)}{\left(\Signal\left[\Sample[\ell_1]\right]-\Signal\left[\Sample[\ell_2]\right]\right)^{2}},\forall h\in\mathbb{R}^{d}
\end{equation*}

where $N\left(h\right)=\left\{\left(\ell_1,\ell_2\right)\in\{1,\ldots,\Samplesize\}^2:\Sample[\ell_1]-\Sample[\ell_2]=h\right\}$ and $|N\left(h\right)|$ is the number of distinct pairs in $N\left(h\right)$. 

\begin{property}
$E[N(h)\hat{\Semivariogram}(h)]=\Semivariogram(h)$
\end{property}

When data are irregularly spaced in $\mathbb{R}^{d}$, we can use
\begin{equation} \label{eq:variogram_hat1}
2\hat{\gamma}_{s}\left(h\left(l\right)\right)=\mathrm{average}\lbrace\left(\Signal\left[\Sample[\ell_1]\right]-\Signal\left[\Sample[\ell_2]\right]\right)^{2}:
\|\Sample[\ell_1]-\Sample[\ell_2]\|\in[h\mp \alpha]\rbrace,
\end{equation}
where the tolerance $\alpha$ is a positive number.

Once the estimated variogram is obtained (or \emph{empirical}), a model is fitted to it in order to achieve a valid variogram. At this stage, we are searching for a valid variogram ``closest'' to the empirical one, and typically we look into a subset of valid variograms $P=\lbrace2\gamma:2\gamma\left(\cdot\right)=2\gamma\left(\cdot;\parampop\right);\parampop\in\parampop\rbrace$. The best element of $P$ can be searched through several good-of-fitness criteria, such as maximum likelihood estimator or Least Squares method.

\begin{figure}[H]
\caption{Naive estimation of the variogram. In the first row, average densities (dotted lines) are compared with the population density (continuos line). In the second row, expected values of the variogram (dotted lines) are compared with population variogram (continuos line).}
\centering
\includegraphics[width=0.7\textwidth]{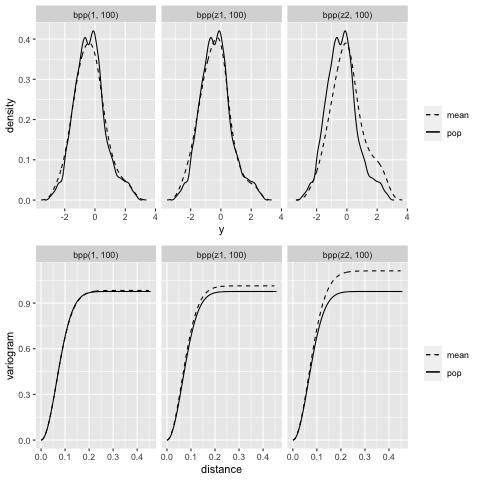}
\label{fig:sim_naive_est}
\end{figure}

In order to investigate the behaviour of the naive estimation, we carried out a Monte Carlo simulation. We simulated one replication of the same isotropic Gaussian process considered in the Example 1, that is $\Signal:\Omega\to(\Pop=[0,1]^2\to\mathbb{R})$, with $\forall \position\in\Pop$, $\mathrm{\mu}[\position]=0$ and with a Gaussian Covariogram with parameters $\parampop_1=5$, $\parampop_2=0.1$. Samples were selected by means of $\mathrm{bpp}(1, 100)$, $\mathrm{bpp}(\desvar_{1}, 100)$, and $\mathrm{bpp}(\desvar_{2}, 100)$, where $\mathbf{z}_{1}=(\log(10)-(0.4^2+0.3^2),0,\sqrt{0.4^2+0.3^2})$ and $\mathbf{z}_{2}=(\log(10)-(0.4^2+0.3^2),0.4,0.3)$, as in the Example 2.
The $\mathrm{bpp}(1, \samplesize)$ is by definition a non-informative sampling mechanism, while the $\mathrm{bpp}(\desvar, \samplesize)$ could introduce bias when the variable $\desvar$ is correlated with the signal $\signal$, as explained in the previous Sections. Note that $\desvar_{1}$ was generated independently from $\signal$ while $\desvar_{2}$ was generated dependently from $\signal$. Indeed, $\mathrm{cor}(\signal, \desvar_{1})\approx -0.08$ and $\mathrm{cor}(\signal, \desvar_{2})\approx 0.75$, where $\mathrm{cor}$ indicates the correlation between the two variables. For each sampling mechanism, $M=1,000$ samples were selected, and variograms estimated by the method of moments and fitted by the Weighted Least Square criterion. These estimated variograms are compared to the variogram obtained by using all the population values, which we call \emph{population variogram}. In the first row of Figure \ref{fig:sim_naive_est}, we compare the averages of the densities obtained by samples selected by means of the three different selection mechanisms. We note that when $\mathrm{bpp}(\desvar_{2}, 100)$ is employed, the average density is different from the \emph{population density}. In the second row of Figure \ref{fig:sim_naive_est}, we compare the expected value of the variograms obtained by the replications to the population variogram. The simulation shows that when the samples are selected by means of $\mathrm{bpp}(1, 100)$ and $\mathrm{bpp}(\desvar_{1}, 100)$, the Monte Carlo expected value of the sample variogram (indicated by the dotted line) is close to the population variogram (indicated by continuos line), while when the samples are selected by means of $\mathrm{bpp}(\desvar_{2}, 100)$, the Monte Carlo expected value of the sample variogram (indicated by the dotted line) is  different from the population variogram (conntinuos line). Therefore, bias is introduced by the selection mechanism.

\subsection{Maximum likelihood estimation}



Here we briefly present a possible solution that take into account the informativeness  of the sample. The proposal is based on the maximum likelihood estimation. In particular, in our setup, the loglikelihood of $\Signal[\position]$ is given by 

\begin{eqnarray} \label{Lnaive}
\likelihood_{\Signal[\position]}\left(\parampop;\signal\right)&=&\log(\density_{\Signal[\position]}(\signal;\parampop))\\&=&
\frac12\log{|\provar_{\position,\position^{'};\parampop}|}
+\frac12(\Signal[\Sample])^{\mathrm{T}}\provar_{\position,\position_{'};\parampop}^{-1}\Signal[\Sample]
\end{eqnarray}

and the naive estimator maximizes \eqref{Lnaive}. Indeed, this maximization ignores the selection mechanism. 

In order to take into account the informativeness of the sample, the \emph{full loglikelihood} should be used, which is composed by three factors: the density ratio, $\densityratio_{\Sample}(\position\mid\signal)$, the distribution of $\Signal$, $\density_{\Signal[\position]}(\signal)$, and the distribution of the sample, $\density_\Sample(\position)$. Therefore the loglikelihood is as in the follow

\begin{eqnarray} \label{fit:ML2}
\likelihood_{\Sample,\Signal[\Sample],\Samplesize}\left(\parampop,\paramnuisance;\position,\signal,\samplesize\right)&=&\log\left(\densityratio_{\Sample}(\position\mid\signal;\parampop,\paramnuisance)\right)+\log(\density_{\Signal[\position]}(\signal;\parampop))+\log\left(\density_\Sample(\position;\parampop,\paramnuisance)\right)
\end{eqnarray}

%% file: sectionC.tex
\section{Conclusions} \label{sec:conclusions}
We have extended the notions of informative selection, population and sample distributions defined by \cite{pfefferman_1992} to a situation where a spatial process is assumed to have generated the whole population. In particular, a Gaussian random field and two point processes that represent the selection mechanism have been considered. An analysis of the naive estimation of the variogram have been performed, and simulation shows that selection bias is introduced when there is dependence between the superpopulation model and the mechanism selection. Finally, we briefly discussed a possible solution in order to correct the naive estimation in presence of informative selection. This solution is based on the use of the maximum likelihood. Future research will be dedicated to a proper development of such proposal.

%% file: Acknowledgements.tex
\section*{Acknowledgements}
We developed an R package that allows to reproduce all the simulation results see \cite{gitSpatialInformativeSelection}. Francesco Pantalone's work was partially funded by the \emph{``International Graduate Research Fellowships''} at Joint Program in Survey Methodology, University of Maryland, College Park, USA.

%% file: sectionA.tex
\newpage
\section{Algebra for the example}
\subsection{Distribution of  $\Desvar$ conditionally on  $\Signal[\position]=\signal$}\label{sec:A.owifjoeij}
Under the condition of the paper example, conditionally on $\Signal[\position]=\signal$, the process $(\alpha+\beta\Signal+\gamma\varepsilon)$ is a Gaussian Process characterised by
$\mathrm{E}\left[(\alpha+\beta\Signal+\gamma\varepsilon)[\position']\right]=\alpha+\beta(\mu+\Sigma_{\position',\position}\Sigma_{\position,\position}^{-1} (\signal-\mu))$ and  $\mathrm{Var}\left[(\alpha+\beta\Signal+\gamma\varepsilon)[\position']\right]=\gamma^2\provar^\varepsilon_{\position',\position'}+\beta^2(\provar_{\position',\position'}-\provar_{\position',\position}\provar_{\position,\position}^{-1}\provar_{\position,\position'})$, and consequently, conditionally on $\Signal[\position]=\signal$, $\Desvar=\exp(\alpha+\beta\Signal+\gamma\varepsilon)$ is a lognormal spatial process.

%
%

\begin{proof}

\begin{eqnarray}\lefteqn{\begin{bmatrix}
\Signal[\position']\\
\Signal[\position]
\end{bmatrix}\sim\mathrm{Normal}\left(\mu,\begin{bmatrix}\provar_{\position',\position'}&\provar_{\position',\position}\\
\provar_{\position,\position'}&
\provar_{\position,\position}
\end{bmatrix}\right)}\\
&\Rightarrow&
\left[
\left.\Signal[\position']
\right|
\Signal[\position]=\signal
\right]
\sim\mathrm{Normal}\left(\mu+\Sigma_{\position',\position}\Sigma_{\position,\position}^{-1} (\signal-\mu),\Sigma_{\position',\position'}-\Sigma_{\position',\position}\Sigma_{\position,\position}^{-1}\Sigma_{\position,\position'}\right)\label{eq:hfqopio}\end{eqnarray}

The independence of $\varepsilon$ and $\Signal$ implies that:
\begin{equation}\label{eq:iushfisduhsd}\left[\left.
\varepsilon[\position']
\right|
\Signal[\position]=\signal
\right]\sim\mathrm{Normal}\left(\mu_\varepsilon,\provar_{\position',\position';\varepsilon}\right),\end{equation} and that the vector obtained by stacking $\Signal[\position]$, $
\Signal[\position']$, and $\varepsilon[\position']$ is normal. 
By combining \eqref{eq:iushfisduhsd} and \eqref{eq:hfqopio} we obtain that:

$$\left[\left.\begin{bmatrix}
\Signal[\position']\\
\varepsilon[\position']
\end{bmatrix}\right|
\Signal[\position]=\signal
\right]\sim\mathrm{Normal}\left(\begin{bmatrix}\mu+\Sigma_{\position',\position}\Sigma_{\position,\position}^{-1} (\signal-\mu)\\0\end{bmatrix},
\begin{bmatrix}\Sigma_{\position',\position'}-\Sigma_{\position',\position}\Sigma_{\position,\position}^{-1}\Sigma_{\position,\position'}&0\\0&\provar_{\position',\position';\varepsilon}\end{bmatrix}\right).$$

We obtain the moments of $Z[\position']=\gamma
\varepsilon[\position']+\alpha+\beta
\Signal[\position']$ conditionally on $\Signal[\position]=\signal$:
$$\mathrm{E}\left[\alpha+\beta
\Signal[\position']+\gamma
\varepsilon[\position']|
\Signal[\position]=\signal\right]
=\alpha+\beta(\mu+\Sigma_{\position',\position}\Sigma_{\position,\position}^{-1} (\signal-\mu)),$$
and 
$$\mathrm{Var}\left[\alpha+\beta
\Signal[\position']+\gamma
\varepsilon\left[\position'\right]|
\Signal[\position]=\signal\right]
=\gamma^2\provar_{\position',\position';\varepsilon}+\beta^2\left(\provar_{\position',\position'}-\provar_{\position',\position}\provar_{\position,\position}^{-1}\provar_{\position,\position'})\right).$$
\end{proof}

\subsection{Algebra for $\provar$}

The denominator of $\rho$ requires the computation of $\int_{U}\Sigma_{\position',\position}~d\dominantU$, which is itself a function of
$\int_{U}\Covariogram(\position'-\position_j)d\dominantU\left(\position'\right)$

For $\position=(\position_1,\position_2)$, such that $\|\position_1-\position_2\|=h$, we have the following:

\begin{eqnarray*}
\Sigma_{\position,\position}&=&
    \begin{bmatrix}\Covariogram(0)&\Covariogram(h)\\\Covariogram(h)&\Covariogram(0)\end{bmatrix}\\
\Sigma_{\position,\position}^{-1}&=&
    (\Covariogram(0)^2-\Covariogram(h)^2)^{-1}\begin{bmatrix}\Covariogram(0)&-\Covariogram(h)\\-\Covariogram(h)&\Covariogram(0)\end{bmatrix}\\
\Sigma_{\position',\position}\Sigma_{\position,\position}^{-1}
    &=&(\Covariogram(0)^2-\Covariogram(h)^2)^{-1}
        \begin{bmatrix}\Covariogram(\position'-\position_1)&\Covariogram(\position'-\position_2)\end{bmatrix}\begin{bmatrix}\Covariogram(0)&-\Covariogram(h)\\-\Covariogram(h)&\Covariogram(0)\end{bmatrix}\\
    &=&\frac{\begin{bmatrix}
        \Covariogram(0)\Covariogram(\position'-\position_1)-\Covariogram(h)\Covariogram(\position'-\position_2)&              \Covariogram(0)\Covariogram(\position'-\position_2)-\Covariogram(h)\Covariogram(\position'-\position_1)\end{bmatrix}}{\Covariogram(0)^2-\Covariogram(h)^2}\\
\Sigma_{\position',\position}&=&\begin{bmatrix}\Covariogram(\position'-\position_1)&\Covariogram(\position'-\position_2)\end{bmatrix}
\end{eqnarray*}